\RequirePackage{ifpdf}
\ifpdf 
\documentclass[pdftex]{sigma}
\else
\documentclass{sigma}
\fi

\newcommand{\Zint}{\mathbb {Z}}

\newcommand{\Rea}{\mathbb {R}}      
\newcommand{\Cplx}{\mathbb {C}}     

\newcommand{\Ha}{\mathfrak {h}}

\begin{document}

\numberwithin{equation}{section}

\allowdisplaybreaks

\renewcommand{\PaperNumber}{038}

\FirstPageHeading

\renewcommand{\thefootnote}{$\star$}

\ShortArticleName{Towards Finite-Gap Integration of the Inozemtsev
Model}

\ArticleName{Towards Finite-Gap Integration\\ of the Inozemtsev
Model\footnote{This paper is a contribution to the Vadim Kuznetsov
Memorial Issue ``Integrable Systems and Related Topics''. The full
collection is available at
\href{http://www.emis.de/journals/SIGMA/kuznetsov.html}{http://www.emis.de/journals/SIGMA/kuznetsov.html}}}

\Author{Kouichi TAKEMURA} 
\AuthorNameForHeading{K. Takemura}

\Address{Department of Mathematical Sciences, Yokohama City
University,\\ 22-2 Seto, Kanazawa-ku, Yokohama 236-0027, Japan}

\Email{\href{mailto:takemura@yokohama-cu.ac.jp}{takemura@yokohama-cu.ac.jp}}

\ArticleDates{Received October 31, 2006, in f\/inal form February
07, 2007; Published online March 02, 2007}

\Abstract{The Inozemtsev model is considered to be a multivaluable
generalization of Heun's equation. We review results on Heun's
equation, the elliptic Calogero--Moser--Suther\-land model and the
Inozemtsev model, and discuss some approaches to the f\/inite-gap
integration for multivariable models.}

\Keywords{f\/inite-gap integration; Inozemtsev model; Heun's
equation; Darboux transformation}

\Classification{81R12; 33E10; 34M35}

\section{Introduction}

Dif\/ferential equations def\/ined on a complex domain frequently
appear in mathematics and physics. One of the most important
dif\/ferential equations is the Gauss hypergeometric
dif\/feren\-tial equation. Mathematically, it is a standard form
of the second-order dif\/ferential equation with three regular
singularities on the Riemann sphere. Global properties of the
solutions, i.e., the monodromy, often play decisive roles in
applications to physics and mathematics.

A canonical form of a Fuchsian equation with four singularities is
given by Heun's equation which is written as
\begin{gather}
\left( \left(\frac{d}{dw}\right) ^2 + \left(
\frac{\gamma}{w}+\frac{\delta }{w-1}+\frac{\epsilon}{w-t}\right)
\frac{d}{dw} +\frac{\alpha \beta w -q}{w(w-1)(w-t)}
\right)\tilde{f}(w)=0, \label{Heun}
\end{gather}
with the condition $\gamma +\delta +\epsilon =\alpha +\beta +1$.
This equation appears in several topics in physics, i.e.,
astrophysics, crystalline materials and so on (see \cite{SL} and
references therein). Although the problem of describing the
monodromy of Heun's equation is much more dif\/f\/icult than that
of the hypergeometric equation, Heun's equation has been studied
from several viewpoints. A~method of f\/inite-gap integration is
available on a study of Heun's equation, and consequently we have
some formulae on the monodromy.

If there exists an odd-order dif\/ferential operator $A= \left(
d/dx \right)^{2g+1} +   \sum\limits_{j=0}^{2g-1}\! $ $ b_j(x)
\left( d/dx \right)^{2g-1-j} $ such that $[A, -d^2/dx^2+q(x)]=0$,
then $q(x)$ is called the algebro-geometric f\/inite-gap
potential. Note that the equation  $[A, -d^2/dx^2+q(x)]=0$ is
equivalent to the function $q(x)$ being a solution to a stationary
higher-order KdV equation (see \cite{DMN}). In our setting,
f\/inite-gap integration is a~method for analysis of the operator
$-d^2/dx^2+q(x)$ where $q(x)$ is an algebro-geometric f\/inite-gap
potential. Originally, the f\/inite-gap property is a notion
related to spectra. Let $H$ be the operator $-d^2/dx^2+q(x)$, and
the set $\sigma _b(H)$ be def\/ined as follows:
\begin{gather*}
E \in \sigma _b(H) \; \Leftrightarrow \mbox{ Every solution to
}(H-E)f(x)=0 \mbox{ is bounded on }x \in \Rea .
\end{gather*}
If the closure of the set $\sigma _b(H)$ can be written as
\begin{gather*}
\overline{\sigma _b(H)}= [E_{0},E_{1}]\cup [E_{2},E_{3}]  \cup
\dots \cup [E_{2g}, \infty ),
\end{gather*}
where $E_0<E_{1}<\cdots <E_{2g}$, i.e., the number of bounded
bands is f\/inite, then $q(x)$ is called the f\/inite-gap
potential. It was established in the 1970s that, under the
assumption that $q(x)$ is a periodic, smooth, real function, the
potential $q(x)$ is f\/inite-gap if and only if $q(x)$ is
algebro-geometric f\/inite-gap.

On the approach by the f\/inite-gap integration for Heun's
equation, it is essential to transform Heun's equation into a form
with the elliptic function. The transformed equation is written as
\begin{gather}
(H-E) f(x)= \left( -\frac{d^2}{dx^2} + \sum_{i=0}^3 l_i(l_i+1)\wp
(x+\omega_i)-E\right) f(x)=0, \label{InoEF0}
\end{gather}
where $\wp (x)$ is the Weierstrass $\wp$-function with periods
$(2\omega _1,2\omega _3)$, $\omega _0(=0)$, $\omega_1$,
$\omega_2(=-\omega _1 -\omega _3)$, $\omega_3$ are half-periods,
and $l_i$ $(i=0,1,2,3)$ are coupling constants. Here the variables
$w$ and $x$ in equations (\ref{Heun}), (\ref{InoEF0}) are related
by $w=(\wp (x)-\wp (\omega _1))/(\wp (\omega _2)-\wp (\omega
_1))$. For details of the transformation, see \cite{Ron,Smi,Tak2}.
The expression in terms of the elliptic function was already
discovered in the 19th century, and some results which may relate
to f\/inite-gap integration were found in that era. For the case
when three of $l_0$, $l_1$, $l_2$, $l_3$ are equal to zero,
equation (\ref{InoEF0}) is called Lam\'e's equation. Ince \cite{I}
established in 1940 that if $n \in \Zint _{\geq 1}$, $\omega _1
\in \Rea \setminus \{ 0 \}$ and $\omega _3 \in \sqrt{-1} \Rea
\setminus \{ 0 \}$, then the potential of Lam\'e's operator
\begin{gather*}
-\frac{d^2}{dx^2}+n(n+1)\wp (x+\omega _3),
\end{gather*}
is f\/inite-gap. In the late 1980s, Treibich and Verdier \cite{TV}
found that the method of f\/inite-gap integration is applicable
for the case $l_0, l_1 , l_2 , l_3 \in \Zint $. Namely, they
showed that the potential in equation (\ref{InoEF0}) is an
algebro-geometric f\/inite-gap potential if $l_i \in \Zint $ for
all $i \in \{ 0,1,2,3 \}$. Therefore the potential $
\sum\limits_{i=0}^3 l_i(l_i+1)\wp (x+\omega_i)$ is called the
Treibich--Verdier potential. Subsequently several others
\cite{GW,Smi,Tak1,Tak3,Tak4,Tak5} have produced more precise
statements and concerned results on this subject. Namely, integral
representations of solutions \cite{Smi,Tak1}, the Bethe Ansatz
\cite{GW,Tak1}, the global monodromy in terms of the hyperelliptic
integrals \cite{Tak3}, the Hermite--Krichever Ansatz~\cite{Tak4}
and a relationship with the Darboux transformation \cite{Tak5}
were studied.

In this paper, we discuss some approaches to f\/inite-gap
integration for multivariable cases. A~multivariable
generalization of Heun's equation is given by the Inozemtsev
model, which is a~generalization of the
Calogero--Moser--Sutherland model. The Inozemtsev model of
ty\-pe~$BC_N$~\cite{Ino} is a quantum mechanical system with
$N$-particles whose Hamiltonian is given~by
\begin{gather*}
 H=-\sum_{j=1}^N\frac{\partial ^2}{\partial x_j^2}+2l(l+1)
\sum_{1\leq j<k\leq N} \left( \wp (x_j-x_k) +\wp (x_j +x_k) \right) \\ 
 \phantom{H=}{}  + \sum_{j=1}^N \sum _{i=0}^3 l_i(l_i+1) \wp(x_j +\omega_i), \nonumber
\end{gather*}
where $l$ and $l_i$ $(i=0,1,2,3)$ are coupling constants. It is
known that the Inozemtsev model of type $BC_N$ is completely
integrable. More precisely, there exist operators of the form
$H_k= \sum\limits_{j=1}^N \left( \frac{\partial }{\partial x_j}
\right) ^{2k} + \mbox{(lower terms)}$ $(k=2,\dots ,N)$ such that
$[H, H_k]=0$ and $[H_{k_1}, H_{k_2}]=0$ $(k, k_1, k_2=2,\dots
,N)$. Note that the Inozemtsev model of type $BC_N$ is a universal
completely integrable model of quantum mechanics with $B_N$
symmetry, which follows from the classif\/ication due to Ochiai,
Oshima and Sekiguchi \cite{OOS}. For the case $N=1$, the potential
coincides with the Treibich--Verdier potential and the spectral
problem for the Inozemtsev model of type $BC_1$ is equivalent to
solving Heun's equation.

By the trigonometric limit $(\tau (=\omega _3/\omega _1)
\rightarrow \sqrt{-1} \infty )$, we obtain the trigonometric
Calogero--Moser--Sutherland model. The trigonometric model is
well-studied by multivariable orthogonal polynomials (i.e., the
Jack polynomial and the multivariable Jacobi polynomial). Vadim
Kuznetsov and his collaborators studied multivariable orthogonal
polynomials from the aspects of separation of variables
\cite{KNS}, the Pfaf\/f lattice \cite{AKM} and the $Q$ operator
\cite{KMS}. Note that in the paper~\cite{KS} relationships among
separation of variables, integral transformations and Lam\'e's
(Heun's) dif\/ferential equation were discussed. Applications of
these technique for models with elliptic potentials are
anticipated.

Although the Inozemtsev model of type $BC_N$ is much more
dif\/f\/icult than the trigonometric one, some approaches
(perturbation from the trigonometric limit, quasi-solvability
etc.) were introduced. Now we hope to develop analysis of this
model by f\/inite-gap integration. Although we can regard several
works to be on this direction, in my opinion, they are still far
from complete understanding of the model. We may consider the
multivariable Darboux transformation as a~possible approach, but
we should develop it in the future.

This paper is organized as follows. In Section \ref{sec:Heun}, we
review the f\/inite-gap integration of Heun's equation. An
approach by the Darboux transformation is introduced. In Section
\ref{sec:resultsICMS}, we collect results on the
Calogero--Moser--Sutherland model and the Inozemtsev model. In
Section \ref{sec:fingapICMS}, we discuss some approaches to
f\/inite-gap integration for those models.

\section{Finite-gap integration of Heun's equation} \label{sec:Heun}

\subsection{Darboux transformations and Heun's equation}

We consider the f\/inite-gap property of Heun's equation in the
elliptic form. It is known that the Treibich--Verdier potential is
algebro-geometric f\/inite-gap, i.e., there exists a
dif\/ferential opera\-tor~$A$ of odd order which commutes with the
operator $H^{(l_0,l_1,l_2,l_3)}$ where
\begin{gather}
H^{(l_0,l_1,l_2,l_3)}= -\frac{d^2}{dx^2} + \sum_{i=0}^3
l_i(l_i+1)\wp (x+\omega_i). \label{Ino}
\end{gather}
In this subsection, we construct an odd-order dif\/ferential
operator $A$ by composing the Darboux--Crum transformation which
we will explain below.

We review the Darboux transformation. Let $\phi _0 (x)$ be an
eigenfunction of the operator $H= -d^2/dx^2 + q(x)$ corresponding
to an eigenvalue $E_0$, i.e.
\begin{gather*}
\left( -\frac{d^2}{dx^2} + q(x) \right) \phi _0(x) = E_0 \phi
_0(x).
\end{gather*}
For this case, the potential $q(x)$ is written as $q(x)=(\phi
_0'(x)/\phi _0(x))'+(\phi _0'(x)/\phi _0(x))^2+E_0$. By setting
$L= d/dx - \phi _0'(x)/\phi _0(x)$ and $L^{\dag }= - d/dx - \phi
_0'(x)/\phi _0(x)$, we have the factorization $H-E_0 =L^{\dag }L$.
We set $\tilde {H}= L L^{\dag }+E_0$. Then we have $\tilde {H}=
-d^2/dx^2 +q(x) - 2(\phi _0'(x)/\phi _0(x))'$ and the relation
$\tilde{H}L=LH$. Hence, if $\phi (x)$ is an eigenfunction of the
operator $H$ corresponding to the eigenvalue $E$, then $L \phi
(x)$ is an eigenfunction of the operator $\tilde{H}$ corresponding
to the eigenvalue $E$. This transformation is called the Darboux
transformation. We generalize the operator $L$ to be the
dif\/ferential operator of higher order in the following
proposition.
\begin{proposition}[cf.~\cite{AST}]\label{prop:GDT}
Suppose that the operator $H=-d^2/dx^2 +q(x)$ preserves an
$n$-di\-men\-sional space $U$ of functions. Let $L$ be the monic
differential operator of order $n$ which annihilates all functions
in $U$, and write
\begin{gather*}
L= \left( \frac{d}{dx} \right) ^n +\sum _{i=1}^n c_i(x) \left(
\frac{d}{dx} \right) ^{n-i}.
\end{gather*}
Set
\begin{gather*}
\tilde{H} = -d^2/dx^2 +q(x)+2c'_1(x).
\end{gather*}
Then we have
\begin{gather*}
\tilde{H}L=LH.
\end{gather*}
\end{proposition}
We call the operator $L$ in Proposition \ref{prop:GDT} the
generalized Darboux transformation or the Darboux--Crum
transformation. For the case $n=1$, let $\phi _0 (x)$ be a
non-zero function in $U$, then $U= \Cplx \phi _0 (x)$, the
operator which annihilates $\phi _0 (x)$ is given by $L=d/dx -
\phi' _0 (x)/\phi _0 (x)$ and the operator $\tilde {H}$ is given
by $\tilde {H}=  H - 2(\phi _0'(x)/\phi _0(x))'$. Hence the
proposition reproduces the Darboux transformation for the case
$n=1$.

We apply Proposition \ref{prop:GDT} to Heun's equation. For this
purpose, we recall the quasi-solvability of Heun's equation. If a
f\/inite-dimensional space is invariant under an action of the
Hamiltonian~$H$ of a model, the model is partially solvable by
diagonalizing the matrix representation of $H$. This situation is
called quasi-solvable. On the quasi-solvability of Heun's
equation, we have
\begin{proposition}[\cite{Tak2}, Proposition 5.1] \label{findim}
Let $\alpha _i$ be a number such that $\alpha _i= -l_i$ or $\alpha
_i= l_i+1$ for all $i\in \{ 0,1,2,3\} $, and set
$d=-\sum\limits_{i=0}^3 \alpha _i /2$. Suppose that $d\in
\Zint_{\geq 0}$, and let $V_{\alpha _0, \alpha _1, \alpha _2,
\alpha _3}$ be the $d+1$-dimensional space spanned by
\begin{gather*}
\left\{ \widehat{\Phi}(\wp (x)) \wp(x)^n\right\} _{n=0, \dots ,d}, 
\end{gather*}
where $\widehat{\Phi}(z)=(z-e_1)^{\alpha _1/2}(z-e_2)^{\alpha
_2/2}(z-e_3)^{\alpha _3/2}$. Then the operator $H^{(l_0,l_1,l_2,
l_3)}$ (see equation~\eqref{Ino}) preserves the space $V_{\alpha
_0, \alpha _1, \alpha _2, \alpha _3}$.
\end{proposition}
By applying Propositions \ref{prop:GDT} and \ref{findim}, we
obtain the following proposition after some calculations:
\begin{proposition}[\cite{Tak5}] \label{thm:HL0123L0123H}
Let $\alpha_i$ be a number such that $\alpha_i= -l_i$ or
$\alpha_i= l_i+1$ for all $i\in \{ 0,1,2,3\} $, and set
$d=-\sum\limits_{i=0}^3 \alpha_i /2$. Suppose that $d\in
\Zint_{\geq 0}$, and let $L_{\alpha _0, \alpha _1, \alpha _2,
\alpha _3}$ be the monic differential operator of order $d+1$
which annihilates the space $V_{\alpha _0, \alpha _1, \alpha _2,
\alpha _3}$. Then we have
\begin{gather*}
H^{(\alpha _0 +d,\alpha _1 +d,\alpha _2 +d,\alpha _3 +d)}  L_{\alpha _0, \alpha _1, \alpha _2, \alpha _3}=L _{\alpha _0, \alpha _1, \alpha _2, \alpha _3} H^{(l_0,l_1,l_2,l_3)}. 
\end{gather*}
\end{proposition}

We construct an odd-degree dif\/ferential operator $A$ which
commutes with $H$ by composing the operators $L_{\alpha _0, \alpha
_1, \alpha _2, \alpha _3}$. We def\/ine the operator
$\tilde{L}_{\alpha _0, \alpha _1, \alpha _2, \alpha _3} $ as
follows:
\begin{gather*}
 \tilde{L}_{\alpha _0, \alpha _1, \alpha _2, \alpha _3}=
\left\{
\begin{array}{ll}
L_{\alpha _0, \alpha _1, \alpha _2, \alpha _3}, & \sum\limits_{i=0}^3 \alpha _i/2 \in \Zint_{\leq 0} ;\vspace{1mm}\\
L_{1-\alpha _0, 1-\alpha _1, 1-\alpha _2, 1-\alpha _3} ,& \sum\limits_{i=0}^3 \alpha _i /2\in \Zint_{\geq 2} ;\vspace{1mm}\\
1 ,& \mbox{otherwise}.
\end{array}
\right. 
\end{gather*}
Set
\begin{alignat*}{3}
& l_0^e= (-l_0+l_1+l_2+l_3)/2, \qquad && l_1^e= (l_0-l_1+l_2+l_3)/2, & \\ 
& l_2^e= (l_0+l_1-l_2+l_3)/2, \qquad && l_3^e= (l_0+l_1+l_2-l_3)/2, & \nonumber \\
& l_0^o= (l_0+l_1+l_2+l_3+1)/2, \qquad && l_1^o= (l_0+l_1-l_2-l_3-1)/2,& \nonumber \\
& l_2^o= (l_0-l_1+l_2-l_3-1)/2, \qquad && l_3^o=
(l_0-l_1-l_2+l_3-1)/2. & \nonumber
\end{alignat*}
The following proposition is proved by applying Proposition
\ref{thm:HL0123L0123H} four times:
\begin{proposition}[\cite{Tak5}] \label{prop:tildeA}
Assume $l_i \in \Zint_{\geq 0}$ for $i=0,1,2,3$. If
$l_0+l_1+l_2+l_3$ is even, we set
\begin{gather}
A= \tilde{L}_{-l_3^e,l_2^e+1,l_1^e+1,-l_0^e}
\tilde{L}_{-l_1,l_0+1,-l_3,l_2+1}
\tilde{L}_{-l_0^e,-l_1^e,l_2^e+1,l_3^e+1}
\tilde{L}_{-l_0,-l_1,-l_2,-l_3} , \label{eq:Ate}
\end{gather}
while if $l_0+l_1+l_2+l_3$ is odd, we set
\begin{gather}
A= \tilde{L}_{l_2^o+1,-l_3^o,-l_0^o,-l_1^o}
\tilde{L}_{-l_1,-l_0,l_3+1,-l_2}
\tilde{L}_{-l_0^o,l_1^o+1,-l_2^o,-l_3^o}
\tilde{L}_{l_0+1,-l_1,-l_2,-l_3} . \label{eq:Ato}
\end{gather}
We then have that the operator $A$ commutes with
$H^{(l_0,l_1,l_2,l_3)}$, i.e.,
\begin{gather*}
AH^{(l_0,l_1,l_2, l_3)}=H^{(l_0,l_1,l_2, l_3)}A. 
\end{gather*}
\end{proposition}
It is known that, if $l_0,l_1,l_2,l_3  \in \Zint_{\geq 0}$, then
there exist four invariant spaces of the operator
$H^{(l_0,l_1,l_2,l_3)}$, which we consider in
Section~\ref{sec:invsp}, and the four operators in Proposition
\ref{prop:tildeA} are related to the four spaces.

Let $k_i$ be the rearrangement of $l_i$ such that $k_0\geq k_1
\geq k_2 \geq k_3 (\geq 0)$. Set
\begin{gather*}
 g=
\left\{
\begin{array}{ll}
k_0, & l_0+l_1+l_2+l_3: \mbox{ even}, \quad k_0+k_3 \geq k_1+k_2 ;\\
(k_0+k_1+k_2-k_3)/2,& l_0+l_1+l_2+l_3: \mbox{ even}, \quad k_0+k_3 < k_1+k_2 ;\\
k_0, & l_0+l_1+l_2+l_3: \mbox{ odd}, \quad k_0 \geq k_1+k_2 +k_3+1;\\
(k_0+k_1+k_2+k_3+1)/2,& l_0+l_1+l_2+l_3: \mbox{ odd}, \quad k_0 <
k_1+k_2+k_3+1 .
\end{array}
\right.\!\!\!
\end{gather*}
Then $g \in \Zint _{\geq 0}$ and the degree of the operator $A$ is
$2g +1$.

For the case $l_0=2$, $l_1=l_2=l_3=0$, we have $g=2$ and the
operator $A$ is expressed as
\begin{gather*}
 L_{2,-1,-1,0} L_{1,-2,1,0} L_{0,2,-1,-1}  L_{-2,0,0,0}  \\
\qquad{}= \left( \frac{d}{dx} +\frac{1}{2}\frac{\wp '(x)}{\wp (x) -e_1}+\frac{1}{2}\frac{\wp '(x)}{\wp (x) -e_2} \right) \left( \frac{d}{dx} +\frac{\wp '(x)}{\wp (x) -e_1}-\frac{1}{2}\frac{\wp '(x)}{\wp (x) -e_2} \right) \nonumber \\
\qquad{}\times \left( \frac{d}{dx} -\frac{\wp '(x)}{\wp (x) -e_1} +\frac{1}{2}\frac{\wp '(x)}{\wp (x) -e_2}+\frac{1}{2}\frac{\wp '(x)}{\wp (x) -e_3} \right) \nonumber \\
\qquad{}\times  \left( \left( \frac{d}{dx} \right) ^2 -\frac{1}{2}\left( \frac{\wp '(x)}{\wp (x) -e_1} +\frac{\wp '(x)}{\wp (x) -e_2}+ \frac{\wp '(x)}{\wp (x) -e_3} \right) \frac{d}{dx} \right) \nonumber \\
\qquad{} = \left( \frac{d}{dx} \right)^{5} -15\wp (x) \left(
\frac{d}{dx} \right)^{3}  -\frac{45}{2} \wp ' (x) \left(
\frac{d}{dx} \right)^{2} -9\left( 5\wp (x)^2-\frac{3}{4}g_2
\right) \frac{d}{dx}. \nonumber
\end{gather*}

\subsection{Application of finite-gap property}

We investigate Heun's equation in the elliptic form
\begin{gather}
\left( H-E\right) f(x)=0, \qquad H=- \frac{d^2}{dx^2}+u(x), \qquad
u(x)= \sum_{i=0}^3 l_i(l_i+1)\wp (x+\omega_i), \label{InoEF}
\end{gather}
by applying the f\/inite-gap integration, which is based on the
commutativity of $H$ $(=-d^2/dx^2+u(x))$ and an odd-order
dif\/ferential operator $A$.

Since $A$ is a monic dif\/ferential operator of order $2g+1$, it
can be expressed in the form
\begin{gather*}
A= (-1)^g \sum_{j=0}^{g} \left( \tilde{a}_j(x)\frac{d}{dx}
+\tilde{b}_j(x) \right) H ^{g-j} ,
\end{gather*}
where $\tilde{a}_0(x)=1$. We have
\begin{gather*}
 0= [(-1)^g A,H]= \sum_{j=0}^{g} \left[ \tilde{a}_j(x)\frac{d}{dx} +\tilde{b}_j(x), - \frac{d^2}{dx^2} +u(x) \right] H^{g-j} \\
\phantom{0} =  \sum_{j=0}^{g} \left( \tilde{a}_j(x)u'(x)+2\tilde{a}'_j(x)\frac{d^2}{dx^2}+(\tilde{a}_j''(x)+2\tilde{b}'_j (x)) \frac{d}{dx}+\tilde{b}''_j (x)\right)H^{g-j} \nonumber \\
\phantom{0} = \sum_{j=0}^{g} \left(2\tilde{a}'_j(x)(-H+u(x))+(\tilde{a}_j''(x)+2\tilde{b}'_j (x)) \frac{d}{dx}+ \tilde{a}_j(x)u'(x)+\tilde{b}''_j (x)\right)H^{g-j} \nonumber \\
\phantom{0} = \sum_{j=0}^{g} \left((\tilde{a}_j''(x)+2\tilde{b}'_j
(x)) \frac{d}{dx}-2\tilde{a}'_{j+1} (x)+ 2\tilde{a}'_j(x)u(x)+
\tilde{a}_j(x)u'(x)+\tilde{b}''_j (x) \right)H^{g-j} . \nonumber
\end{gather*}
Hence we obtain
\begin{gather*}
 \tilde{b}_j(x)= -\tilde{a}'_j (x)/2 +c_j, \qquad \tilde{a}'''_j(x)-4u(x)\tilde{a}'_j(x)+4\tilde{a}'_{j+1}(x)-2u'(x)\tilde{a}_j(x)=0
\end{gather*}
for some constants $c_j$ $(j=0,\dots ,g)$. Therefore we have
\begin{proposition} \label{prop:tA}
Set $\tilde{a}_0(x)=1$ and $\tilde{a}_{g+1}(x)=0$. The operator
$A$ may be expressed in the form
\begin{gather}
A= (-1)^g\left[ \sum_{j=0}^{g} \left\{
\tilde{a}_j(x)\frac{d}{dx}-\frac{1}{2} \left( \frac{d}{dx}
\tilde{a}_j(x) \right) \right\} H ^{g-j} + \sum_{j=0}^{g} c_j
H^{g-j} \right] , \label{Adef0}
\end{gather}
for some functions $\tilde{a}_j(x)$ $(j=1,\dots ,g)$ and constants
$c_j$ $(j=0,\dots ,g)$, where the functions~$\tilde{a}_j(x)$
$(j=0,\dots ,g)$ satisfy
\begin{gather}
\tilde{a}'''_j(x)-4u(x)\tilde{a}'_j(x)+4\tilde{a}'_{j+1}(x)-2u'(x)\tilde{a}_j(x)=0.
\label{eq:a'''}
\end{gather}
\end{proposition}
We def\/ine a function $\Xi (x,E)$ which plays the important role
for the solutions and the monodromy of Heun's equation. Set
\begin{gather}
\Xi (x,E) = \sum_{i=0}^{g} \tilde{a}_{g-i}(x) E^i. \label{Xiag}
\end{gather}
It follows from equation (\ref{eq:a'''}) that $\Xi (x,E)$
satisf\/ies a dif\/ferential equation satisfied by products of any
pair of the solutions to equation (\ref{InoEF}), i.e.,
\begin{gather}
 \left( \frac{d^3}{dx^3}-4\left( u(x)-E\right)\frac{d}{dx}-2u'(x)  \right) \Xi (x,E)=0.
\label{prodDE}
\end{gather}

On the basis of Proposition \ref{prop:tA} and the function $\Xi
(x,E)$ in equation (\ref{Xiag}), we have
\begin{proposition}[\cite{Tak5}] \label{prop:XixE}
(i) The constants $c_j$ $(j=1,\dots ,g)$ in equation \eqref{Adef0}
are all zero.

(ii) The function $\Xi (x,E)$ is even doubly-periodic and
expressed as
\begin{gather}
\Xi (x,E)=c_0(E)+\sum_{i=0}^3 \sum_{j=0}^{l_i-1} b^{(i)}_j (E)\wp
(x+\omega_i)^{l_i-j}, \label{Fx}
\end{gather}
where the coefficients $c_0(E)$ and $b^{(i)}_j(E)$ are polynomials
in $E$. The coefficients do not have common divisors and the
polynomial $c_0(E)$ is monic. We have $g=\deg_E c_0(E)$ and the
coefficients satisfy $\deg _E b^{(i)}_j(E)<g$ for all $i$ and $j$.
\end{proposition}

For the case $l_0=2$, $l_1=l_2=l_3=0$, we have
\begin{gather*}
 \Xi (x,E)= E^2+3E \wp (x) + 9\wp (x)^2 -\tfrac{9}{4} g_2.
\end{gather*}
where $e_i =\wp (\omega _i )$ and $g_2= -4(e_1e_2+e_2e_3+e_3e_1)$.

Note that the function $\Xi (x,E)$ can be also obtained as the
function satisfying equation~(\ref{prodDE}) and Proposition
\ref{prop:XixE} (ii) (see \cite{Tak1}).

We can derive an integral formula for a solution to
equation~(\ref{InoEF}) in terms of the doubly periodic function
$\Xi(x,E)$. Set
\begin{gather}
Q(E)=  \Xi (x,E)^2\left( E- \sum_{i=0}^3 l_i(l_i+1)\wp (x+\omega_i)\right) \nonumber\\
\phantom{Q(E)=}{} +\frac{1}{2}\Xi (x,E)\frac{d^2\Xi
(x,E)}{dx^2}-\frac{1}{4}\left(\frac{d\Xi (x,E)}{dx} \right)^2.
\label{const}
\end{gather}
It is shown by dif\/ferentiating the right-hand side of equation
(\ref{const}) and applying equation~(\ref{prodDE}) that $Q(E)$ is
independent of $x$. Thus $Q(E)$ is a monic polynomial in $E$ of
degree $2g+1$, which follows from the expression for $\Xi (x,E)$
given by equation (\ref{Fx}). For the case $l_0=2$,
$l_1=l_2=l_3=0$, we have
\begin{gather*}
 Q(E)=(E^2 -3g_2)\prod _{i=1}^3 (E-3e_i).
\end{gather*}
The following proposition on the integral representation of a
solution to equation (\ref{InoEF}) was obtained in \cite{Tak1}:
\begin{proposition}[\cite{Tak1}, Proposition 3.7]
Let $\Xi (x,E)$ be the doubly periodic function defined in
equation \eqref{Xiag} and $Q(E)$ be the monic polynomial defined
in equation \eqref{const}. Then the function
\begin{gather*}
\Lambda (x,E)=\sqrt{\Xi (x,E)}\exp \int \frac{\sqrt{-Q(E)}dx}{\Xi
(x,E)}
\end{gather*}
is a solution to the differential equation \eqref{InoEF}.
\end{proposition}

If the value $E$ satisf\/ies $Q(E) \neq 0$, then the functions
$\Lambda (x,E)$ and $\Lambda (-x,E)$ form a basis of solutions to
equation~(\ref{InoEF}). Since equation~(\ref{InoEF}) is
doubly-periodic, the functions $\Lambda (x +2\omega _k,E)$ and
$\Lambda (-(x+2\omega _k),E)$ are also solutions to
equation~(\ref{InoEF}). We consider the monodromy on the functions
$\Lambda (x,E)$ and $\Lambda (-x,E)$. Note that, if $\Lambda (x
+2\omega _k,E)$ is expressed as $B(E) \Lambda (x,E)$, then
$\Lambda (-(x +2\omega _k),E)$ is expressed as $B(E)^{-1} \Lambda
(-x,E)$. We will express $B(E)$ as a hyperelliptic integral of
second kind. We rewrite the function $\Xi (x,E)$ and def\/ine
$a(E)$ as follows:
\begin{gather*}
\Xi (x,E)=c(E)+\sum_{i=0}^3 \sum_{j=0}^{l_i-1 } a^{(i)}_j
(E)\left( \frac{d}{dx} \right) ^{2j} \wp (x+\omega_i), \qquad
a(E)=\sum _{i=0}^3 a^{(i)} _0 (E).
\end{gather*}
\begin{proposition}[\cite{Tak3}, Theorem 3.7] \label{thm:conj3}
Assume $l_i \in \Zint_{\geq 0}$ ($i=0,1,2,3$). Let $E_0$ be a
value such that $Q(E_0)=0$. Then $\Lambda (x+2\omega
_k,E_0)=(-1)^{q_k} \Lambda (x,E_0)$ for $q_k \in \{0,1\}$
$(k=1,3)$ and we have
\begin{gather*}
\Lambda (x+2\omega _k,E)=(-1)^{q_k} \Lambda (x,E) \exp \left(
-\frac{1}{2} \int_{E_0}^{E}\frac{ -2\eta _k a(\tilde{E}) +2\omega
_k c(\tilde{E}) }{\sqrt{-Q(\tilde{E})}} d\tilde{E}\right) ,
\end{gather*}
where $\eta _k=\zeta (\omega _k)$ $(k=1,3)$ and $\zeta (x)$ is the
Weierstrass zeta function.
\end{proposition}

For the case $l_0=2$, $l_1=l_2=l_3=0$, we set $E_0 =\sqrt{3g_2}$.
Then $q_1 =q_3=0$ and the function~$a(E)$ and~$c(E)$ are
determined as
\begin{gather*}
c(E)=E^2-\tfrac{3}{2}g_2, \qquad a_0 (E) = 3E.
\end{gather*}
Hence we have
\begin{gather*}
\Lambda (x+2\omega _k,E)= \Lambda (x,E) \exp \left( -\frac{1}{2}
\int_{\sqrt{3g_2}}^{E}\frac{ -6 \tilde{E}\eta _k
+(2\tilde{E}^2-3g_2 ) \omega _k  }{\sqrt{-(\tilde{E}^2
-3g_2)\prod\limits_{i=1}^3 (\tilde{E}-3e_i) }} d\tilde{E}\right).
\end{gather*}

We review the propositions related with the Bethe Ansatz
(Proposition \ref{prop:BA}) and the Hermite--Krichever Ansatz
(Proposition~\ref{prop:HKA}), which are also reductions of the
f\/inite-gap property.
\begin{proposition}[\cite{Tak1}, Theorem 3.12] \label{prop:BA}
(i) If the value $E$ satisfies $Q(E)\neq 0$, then there exists
$t_1, \dots, t_n$ and $C$ such that $t_j\neq t_{j'}$ $(j\neq j')$,
$t_j \not \in \omega _1 \Zint + \omega _3 \Zint $ and $\Lambda
(x,E)$ is expressed as
\begin{gather}
 \Lambda (x,E)= C \frac{\prod\limits_{j=1}^l \sigma(x+t_j)}{\sigma(x)^{l_0}\sigma_1(x)^{l_1}\sigma_2(x)^{l_2}\sigma_3(x)^{l_3}}\exp \left(-x\sum_{i=1}^l \zeta(t_j)\right) , \label{BA}
\end{gather}
where $\sigma (x)$ is the Weierstrass sigma function and $\sigma_
i(x)$ $(i=1,2,3)$ are the co-sigma functions defined by
\begin{gather*}
 \sigma_i(z)=\exp (-\eta_i z)\sigma(z+\omega_i)/\sigma(\omega _i).
\end{gather*}

(ii) The function
\begin{gather}
\tilde{\Lambda }(x)=\frac{\prod\limits_{j=1}^l
\sigma(x+t_j)}{\sigma(x)^{l_0}\sigma_1(x)^{l_1}\sigma_2(x)^{l_2}\sigma_3(x)^{l_3}}\exp(cx),
\label{BV}
\end{gather}
with the condition $t_j\neq t_{j'}$ $(j\neq j')$ and $t_j \not \in
\omega _1 \Zint + \omega _3 \Zint $ is an eigenfunction of the
operator $H$ (see equation~\eqref{InoEF}), if and only if  $t_j$
$(j=1,\dots ,l)$ and $c$ satisfy the relations,
\begin{gather}
 \sum_{k\neq j} \zeta(-t_j+t_k)-l_0\zeta(-t_j)-\sum_{i=1}^3l_i(\zeta(-t_j+\omega _i )-\zeta(\omega_i))= -c,\qquad (j=1, \dots ,l),
\label{BAeq}\\
 (1-\delta_{l_0,0})\left(c+\sum_{j=1}^l \zeta(t_j) \right)=0, \nonumber \\
 (1-\delta_{l_i,0})\left(c+l\zeta(\omega_i)+\sum_{j=1}^l \zeta(-\omega_i+t_j)\right)=0, \qquad (i=1,2,3). \nonumber
\end{gather}
The eigenvalue $E$ is given by
\begin{gather*}
 E=-c^2+(l_0l_1+l_2l_3)e_1+(l_0l_2+l_1l_3)e_2+(l_0l_3+l_1l_2)e_3-\sum_{i=1}^3l_i\eta_i (2c+l\eta_i)\\
\phantom{E=}{} -\sum_{j=1}^l\sum_{i=0}^3
l_i(\wp(t_j-\omega_i)-\zeta(t_j-\omega_i)^2)+\sum_{j<k}(\wp(t_j-t_k)-\zeta(t_j-t_k)^2).
\nonumber
\end{gather*}
\end{proposition}

Equation~(\ref{BAeq}) is called the Bethe Ansatz equation for the
Inozemtsev model of type $BC_1$ (see \cite{Tak1}). Note that
Gesztesy and Weikard \cite{GW} obtained similar results. The
monodromy of the function $\tilde{\Lambda }(x)$ in
equation~(\ref{BV}) is written as
\begin{gather*}
 \tilde{\Lambda }(x+2\omega _k) = \exp (2\eta _k (t_1 +\cdots +t_l ) +2  \omega _k (c- \zeta(t_1) -\dots - \zeta(t_l))) \tilde{\Lambda }(x) 
\end{gather*}
for $k=1,2,3$.

In order to describe the proposition on the Hermite--Krichever
Ansatz, we def\/ine
\begin{gather*}
\Phi _i(x,\alpha )= \frac{\sigma (x+\omega _i -\alpha ) }{ \sigma
(x+\omega _i )} \exp (\zeta( \alpha )x), \qquad (i=0,1,2,3).
\end{gather*}
\begin{proposition}[\cite{Tak4}] \label{prop:HKA}
There exist polynomials $P_1 (E) ,\dots ,P_6(E)$ such that, if
$P_2(E) \neq 0$, then $\Lambda (x,E)$ is written as
\begin{gather}
 \Lambda (x,E) = \exp \left( \kappa x \right) \left( \sum _{i=0}^3 \sum_{j=0}^{l_i-1} \tilde{b} ^{(i)}_j \left( \frac{d}{dx} \right) ^{j} \Phi _i(x, \alpha ) \right)
\label{Lalpha}
\end{gather}
for some values $\tilde{b} ^{(i)}_j$ $(i=0,\dots ,3, \: j= 0,\dots
,l_i-1)$, $\alpha $ and $\kappa $. The values $\alpha $ and
$\kappa $ are expressed as
\begin{gather*}
 \wp (\alpha ) =\frac{P_1 (E)}{ P_2 (E)}, \qquad \wp ' (\alpha ) =\frac{P_3 (E)}{P_4 (E)} \sqrt{-Q(E)} , \qquad \kappa  = \frac{P_5 (E)}{P_6 (E)} \sqrt{-Q(E)}.
\end{gather*}
For the periodicity of the function $\Lambda (x,E) $, we have
\begin{gather*}
 \Lambda (x+2\omega _k,E) = \exp (-2\eta _k \alpha +2\omega _k \zeta (\alpha ) +2 \kappa \omega _k ) \Lambda (x,E)  
\end{gather*}
for $k=1,3$.
\end{proposition}
Note that $\alpha $ in equation~(\ref{Lalpha}) and $t_j$ in
equation~(\ref{BA}) satisfy the relation $\alpha = -\sum\limits
_{j=1}^l t_j$. To calculate the polynomials $P_1 (E) ,\dots
,P_6(E)$, it is effective to apply the notions ``twisted Heun
polynomial'' and ``theta-twisted Heun polynomial''
(see~\cite{Tak4}).

If $l_0=2$, $l_1=l_2=l_3=0$, then the values $\alpha $ and $\kappa
$ are expressed as
\begin{gather*}
 \wp( \alpha )= e_1 -\frac{(E-3e_1)(E+6e_1)^2}{9(E^2-3g_2)}, \qquad \kappa =\frac{2}{3(E^2-3g_2)}\sqrt{-Q(E)}.
\end{gather*}

\subsection{Relationship among commuting operators} \label{sec:invsp}

We review a relationship among the operators $H$, $A$, the
polynomial $Q(E)$ and the invariant subspaces. On the operators
$H$ and $A$, we have the following relation:
\begin{proposition}[\cite{Tak3}, Proposition 3.2] \label{prop:algrel}
Let $H$ and $A$ be the operators defined by equation~\eqref{InoEF}
and equations~\eqref{eq:Ate}, \eqref{eq:Ato}, and $Q(E)$ be the
polynomial defined in equation~\eqref{const}. Then
\begin{gather*}
A^2+Q(H)=0. 
\end{gather*}
\end{proposition}
It is known that, if $l_0,l_1,l_2,l_3  \in \Zint_{\geq 0}$, then
there exist four invariant subspaces with respect to the action of
the operator $H$. We describe the spaces more precisely. Let
$V_{\alpha _0, \alpha _1, \alpha _2, \alpha _3}$ be the space
def\/ined in Proposition~\ref{findim} and
\begin{gather*}
 U_{\alpha _0, \alpha _1, \alpha _2, \alpha _3}=
\left\{
\begin{array}{ll}
V_{\alpha _0, \alpha _1, \alpha _2, \alpha _3}, & \sum\limits_{i=0}^3 \alpha _i/2 \in \Zint_{\leq 0} ;\vspace{1mm}\\
V_{1-\alpha _0, 1-\alpha _1, 1-\alpha _2, 1-\alpha _3} ,& \sum\limits_{i=0}^3 \alpha _i /2\in \Zint_{\geq 2} ;\vspace{1mm}\\
\{ 0 \} ,& \mbox{otherwise}.
\end{array}
\right. 
\end{gather*}
If $l_0,l_1,l_2,l_3  \in \Zint_{\geq 0}$ and $l_0 +l_1 +l_2 +l_3$
is even, then the operator $H$ preserves the space
\begin{gather}
V=U_{-l_0,-l_1,-l_2,-l_3}\oplus U_{-l_0 ,-l_1,l_2+1 ,l_3+1}\oplus
U_{-l_0,l_1+1,-l_2,l_3+1}\oplus  U_{-l_0 ,l_1+1,l_2+1,-l_3},
\label{eq:Ve}
\end{gather}
and also preserves the components in equation (\ref{eq:Ve}). If
$l_0,l_1,l_2,l_3  \in \Zint_{\geq 0}$ and $l_0 +l_1 +l_2 +l_3$ is
odd, then the operator $H$ preserves the space
\begin{gather}
V=U_{-l_0,-l_1,-l_2,l_3+1}\oplus U_{-l_0 ,-l_1,l_2+1,-l_3}\oplus
U_{-l_0 ,l_1+1,-l_2,-l_3}\oplus U_{l_0+1,-l_1,-l_2,-l_3},
\label{eq:Vo}
\end{gather}
and also preserves the components in equation~(\ref{eq:Vo}). Then
we have
\begin{proposition}[\cite{Tak5}]
(i) The operator $A$ annihilates any elements in the space $V$.

(ii) The monic characteristic polynomial of the space $V$ with
respect to the action of $H$ coincides with $Q(E)$.
\end{proposition}

Note that the operator $A$ was constructed by composing the
generalized Darboux transformations which are related to the
spaces in the components of equation~(\ref{eq:Ve}) or
equation~(\ref{eq:Vo}).

\section[Results on the Calogero-Moser-Sutherland model and the Inozemtsev model]{Results on the Calogero--Moser--Sutherland model\\ and the Inozemtsev model} \label{sec:resultsICMS}

We are going to consider multidimensional generalizations of
Lam\'e's equation and  Heun's equation in the elliptic form. For
this purpose, we introduce the quantum mechanical systems.

\subsection[The elliptic Calogero-Moser-Sutherland model]{The elliptic Calogero--Moser--Sutherland model} \label{sec:resultAN-1}

The elliptic Calogero--Moser--Sutherland model (or the elliptic
Olshanetsky--Perelomov mo\-del~\cite{OP}) of type $A_{N-1}$ is a
quantum many-body system whose Hamiltonian is given as follows:
\begin{gather*}
H =- \frac{1}{2} \sum_{i=1}^{N} \frac{\partial ^{2}}{\partial
x_{i}^{2}} + l (l+1) \sum_{1 \leq i<j \leq N} \wp ( x_{i}-x_{j}),
\end{gather*}
where $\wp(x)$ is the Weierstrass elliptic function. For the case
$N=2$, the model reproduces Lam\'e's equation by setting
$x_1-x_2=x$ and restricting to the line $x_1+x_2=0$.

This model is known to be completely integrable, i.e., there exist
$N$-algebraically independent commuting operators $P_k$
$(k=1,\dots ,N)$ which commute with the Hamiltonian $H$. Namely,
by setting
\begin{gather}
P_k=\sum _{0\leq j \leq [k/2]}\frac{(l(l+1))^j}{2^j j! (k-2j)!} \sum_{\sigma \in S_N} \sigma (\wp (x_1-x_2) \wp (x_3-x_4) \cdots\nonumber\\
\phantom{P_k=}{}\times \wp (x_{2j-1}-x_{2j}) \partial _{2j+1
}\partial _{2j+2 } \cdots \partial _k ), \label{eq:AN-1Pk}
\end{gather}
where $S_N$ is the symmetric group, $[x]$ is the integral part of
$x$ and $\partial _{i } = \partial /  \partial x_i$, we have
$[P_k, H]= 0$ $(1\leq k\leq N)$ and $[P_k, P_{k'}]= 0$  $(1\leq k,
k'\leq N)$ (see \cite{OOS}). The Hamiltonian $H$ is expressed as
$H=P_2-P_1^2/2$.

By the trigonometric limit $(\tau \rightarrow \sqrt{-1} \infty)$
of the elliptic Calogero--Moser--Sutherland model where $(1,\tau)$
is the basic periods of the elliptic function, we obtain (up to an
additive scalar) the Hamiltonian of the trigonometric
Calogero--Moser--Sutherland model,
\begin{gather*}
 H_{\rm trig} =- \frac{1}{2} \sum_{i=1}^{N}
\frac{\partial ^{2}}{\partial x_{i}^{2}} +\pi ^{2} l(l+1) \sum_{1
\leq i<j \leq N} \frac{1}{\sin ^{2} (\pi  (x_{i}-x_{j}))}.
\end{gather*}

The eigenstates of the Calogero--Sutherland model are described by
the Jack polynomial $J_{\lambda}^{(\frac{1}{l+1})}(X)$ ($\lambda
\in {\cal{M}}_N$) (see \cite{Sta}), i.e.,
\begin{gather*}
 H_{\rm trig} ( J_{\lambda}^{(\frac{1}{l+1})}(X)\Delta (X)^{l+1} ) = (e_0+ 2\pi^2 E_{\lambda }^{[\frac{1}{l+1}]})J_{\lambda}^{(\frac{1}{l+1})}(X)\Delta (X)^{l+1},
\end{gather*}
where $X_{i}=\exp \left( 2 \pi \sqrt{-1} x_{i} \right)$,
${\cal{M}}_N=\{ \lambda = (\lambda _{1}  ,  \lambda _{2} , \dots ,
\lambda _{N}   )| i>j \Rightarrow \lambda _{i} - \lambda_j \in
\Zint_{\geq 0} \}$, $e_0= \frac{1}{6} \pi ^{2} (l+1)^{2}
N(N^{2}-1)$, $\Delta (X)=(X_{1}X_{2} \cdots X_{N}
)^{\frac{1-N}{2}} \prod \limits_{i<j} (X_{i}-X_{j})$ and
$E_{\lambda }^{[\alpha ]}= \sum\limits_{i=1}^{N} \lambda _i^2 +
\sum\limits_{i=1}^{N} \frac{N+1-2i}{\alpha } \lambda _i $. In
particular, the ground-state is given by $\Delta (X)^{l+1}$.
Several properties of the Jack polynomial were studied. Vadim
Kuznetsov and his collaborators studied the Jack polynomial and
related polynomials from the aspects of separation of variable
\cite{KNS}, the Pfaf\/f lattice \cite{AKM} and the $Q$
operator~\cite{KMS}.

In contrast with the trigonometric models, the elliptic models are
less investigated and the spectra or the eigenfunctions are not
suf\/f\/iciently analyzed. There is, however, some important
progress due to Felder and Varchenko. They introduced the Bethe
Ansatz method for the $N$-par\-ticle elliptic Calogero--Moser
model with the coupling constant $l$ a positive integer. Note that
Hermite essentially introduced the Bethe Ansatz method for the
case $N=2$ and $l\in \Zint $ (see \cite{WW}), and Dittrich and
Inozemtsev \cite{DI} did it for the case $N=3$ and $l=1$ in a
dif\/ferent representation.

Fix the parameters $N$ and $l$. We set $m=lN(N-1)/2$. Let $c:
\{1,\dots ,m\} \rightarrow \{1, \dots ,N\}$ be the unique
non-decreasing function such that $c^{-1}(j)$ has $(N-j)l$
elements. Let $\epsilon _i$ $(1\leq i \leq N)$ be an orthonormal
basis of $\Rea ^{N}$ with an inner product $(\cdot , \cdot )$. Set
$ \alpha _i =\epsilon_i - \epsilon_{i+1}$, $\Ha ^{*}= \{
\sum\limits_{i=1}^{N}x _i \epsilon_i | \sum\limits_{i=1}^{N}x _i=0
\}$, $p_i=i(2N-i-1)l/2$ and $V_{i}= \{ p_{i-1}+1, p_{i-1}+2, \dots
, p_i \}$ $(1\leq i\leq N-1)$. Let $W$ be the set of maps $w=(w_1,
\dots ,w_N)$ $(w_i: \; V_i \rightarrow \{i,i+1,\dots ,N-1 \} )$
such that $\# \{ w_i^{-1}(j)\} =l$ for $1\leq i\leq j\leq N-1$.
For $w=(w_1, \dots ,w_{N-1}) \in W$, let $F_{w}$ be the set of
maps $f=(f_1, \dots ,f_{N-2})$ $(f_{i}: \; V_{i+1} \rightarrow
V_i)$ such that (i) $f_{i}$ is injective (ii) If $f_{i}(x)=y$ then
$w_{i+1}(x)=w_i(y)$. Set
\begin{gather*}
 \theta_1(x)=2\sum_{n=1}^{\infty} (-1)^{n-1} \exp \big(\tau \pi \sqrt{-1}(n-1/2)^2\big) \sin(2n-1)\pi x,
 \\
 \theta (x)=\frac{\theta_1(x)}{\theta_1 '(0)}, \qquad \sigma_{\lambda }(x)=\frac{\theta '(0)\theta(x-\lambda)}{\theta(x)\theta(\lambda)}. \nonumber
\end{gather*}

For $\xi \in \Ha ^*$, we introduce the functions $\Phi
_{\tau}(t_1, \dots , t_m)$ and $\omega (t;x)$ as follows
\begin{gather*}
 \Phi _{\tau}(t_1, \dots , t_m)=
 e^{2\pi \sqrt{-1} (\xi,\sum_j t_j\alpha_{c(j)}) } \\
\phantom{\Phi _{\tau}(t_1, \dots , t_m)=}{} \times\prod_{1\leq j
\leq (N-1)l} \theta(t_j) ^{-lN} \prod_{c(i)=c(j)\atop{i<
j}}\theta(t_i-t_j )^{2}
\prod_{|c(i)-c(j)|=1\atop{i< j}}\theta(t_i-t_j )^{-1}, \nonumber \\
 \omega (t;x)=e^{2\pi \sqrt{-1} (\xi ,\sum_{i} x_i \epsilon_i)}
\sum_{w \in W} \sum_{f \in F_w} \prod_{i=1}^{N-1}  \prod_{k=
p_{i-1}+1}^{p_i} \sigma_{x_i -x_{w_i(k)+1}}(t_k-t_{f_i(k)}) ,
\nonumber
\end{gather*}
where $t_0=0$, $f_0(k)=0$. Then we have
\begin{proposition}[\cite{FVKZB, FVthr, FRV}] \label{prop:BAAN-1}
If $(t^0_1, \dots , t^0_m)$ satisfy the following Bethe Ansatz
equations,
\begin{gather*}
\frac{\partial \Phi_{\tau}}{\partial t_i}| _{(t^0_1, \dots , t^0_m)}=0\qquad  (1\leq i \leq m), 
\end{gather*}
the function $\omega (t^0;x)$ is an eigenfunction of the
Hamiltonian $H$ with the eigenvalue
\begin{gather*}
2\pi ^2(\xi,\xi) -2\pi \sqrt{-1} \frac{\partial }{\partial \tau}
S(t_{1}^0 ,\dots,t_{m}^0 ; \tau ) -l(l+1)(N-1)N\eta ,
\end{gather*}
where
\begin{gather*}
S(t_1,\dots,t_m; \tau )= \sum\limits_{i<j}(\alpha_{c(i)}, \alpha_{c(j)}) \log \theta (t_i -t_j)- \sum\limits_{c(i)=1} lN\log \theta (t_i),\\
\eta = \pi^2 \left(\frac{1}{6}-4\sum\limits_{n=1}^{\infty} \frac{p
^{n}}{1-p ^{n}}\right) \qquad {\rm and}\qquad p=\exp \big(2\pi
\sqrt{-1} \tau \big).
\end{gather*}
\end{proposition}
Therefore, if we f\/ind solutions to the Bethe Ansatz equations,
we can investigate the Calogero--Moser--Sutherland model in more
detail. There are two things to be considered for applying
Proposition \ref{prop:BAAN-1} to the spectral problem of the
elliptic Calogero--Moser--Sutherland model. The f\/irst one is to
f\/ind the condition when the eigenfunctions obtained by the Bethe
Ansatz method are connected to square-integrable eigenstates and
the second one is how the solutions of the Bethe Ansatz equation
behave.

On the f\/irst question, the condition is described as the
parameter $\xi $ belonging to some lattice (the weight lattice of
type $A_{N-1}$). By symmetrizing or anti-symmetrizing the function
$\omega (t^0;x)$, we obtain square-integrable eigenstates,
although we must check that they are identically zero or not.

On the second question, we consider the solution at $p= \exp (2\pi
\sqrt{-1} \tau )=0$ (the case of the trigonometric limit $\tau
\rightarrow \sqrt{-1} \infty $) and look into the behavior where
$p$ is near $0$, because it is hopeful to solve the Bethe Ansatz
equations for the trigonometric case in contrast to being hopeless
directly for the elliptic case. A key tool to connect the
trigonomertic solutions to the elliptic solutions is the implicit
function theorem. Thus we construct the square-integrable
eigenstates and obtain the main result in \cite{Tak}, which gives
a suf\/f\/icient condition for regular convergence of the
perturbation expansion. In particular, for the case $N=2$, $l \in
\Zint _{\geq 1}$ and the case $N=3$, $l=1$, we have convergence of
the perturbation series for all eigenstates related to the Jack
polynomial.

Note that this idea can be interpreted to consider the elliptic
Calogero--Moser--Sutherland model by perturbation from the
trigonometric Calogero--Moser--Sutherland model. Convergence for
the general cases was proved in \cite{KT} by another method.
Namely, by applying Kato--Rellich theory, we have convergence of
the perturbation series in $p$ for $l\geq 0$ and arbitrary $N$.
The eigenvalues and the eigenfunctions are calculated as power
series by a standard algorithm of perturbation. Remark that
Fernandez, Garcia and Perelomov \cite{FGP} derived a fully
explicit formula for second order in $p$, and Langmann
\cite{Lang,Lang06} obtained another algorithm for constructing the
eigenfunctions and the eigenvalues as formal power series of $p$,
which also gives a formula for all orders in $p$.

On the Bethe Ansatz for the elliptic Calogero--Moser--Sutherland
model, there are some problems to be solved. For example, it has
not been shown at the moment of writing that the eigenfunction
$\omega (t^0,x)$ written in the form of the Bethe Ansatz is also
an eigenfunction of the higher commuting operators $P_3, \dots
,P_N$ (see also \cite{R}).

\subsection{The Inozemtsev model}

We now introduce a quantum mechanical system that is a
multidimensional generalization of Heun's equation in the elliptic
form.

The Inozemtsev model of type $BC_N$ \cite{Ino} is a quantum
mechanical system with $N$-particles whose Hamiltonian is given by
\begin{gather}
 H=-\sum_{j=1}^N\frac{\partial ^2}{\partial x_j^2}+2l(l+1)\sum_{1\leq j<k\leq N} \left( \wp (x_j-x_k) +\wp (x_j +x_k) \right) \label{InoHam} \\
\phantom{H=}{} + \sum_{j=1}^N \sum _{i=0}^3 l_i(l_i+1) \wp(x_j
+\omega_i), \nonumber
\end{gather}
where $l$ and $l_i$ $(i=0,1,2,3)$ are coupling constants. This is
also a generalization of the elliptic Calogero--Moser--Sutherland
model of type $BC_N$.

The Inozemtsev model of type $BC_N$ is completely integrable,
 i.e., there exist $N$ algebraically independent mutually commuting dif\/ferential operators $P_k$ $(k=1, \dots , N)$ (higher commuting Hamiltonians) which commute with the Hamiltonian of the model,
and Oshima \cite{O} described the commuting operators explicitly.
Note that the Inozemtsev model of type $BC_N$ (resp. the elliptic
Calogero--Moser--Sutherland model of type $A_N$) is a universal
completely integrable model of quantum mechanics with the symmetry
of the Weyl group of type $B_N$ (resp. type $A_N$), which follows
from the classif\/ication due to Ochiai, Oshima and Sekiguchi
\cite{OOS,OS}. For the case $N=1$, the operator (\ref{InoHam})
appears in the elliptic form of Heun's equation (\ref{InoEF}).
Therefore the Inozemtsev model of type $BC_N$ is regarded as a
multidimensional generalization of Heun's equation.

On the trigonometric limit $\tau \rightarrow \sqrt{-1} \infty $,
we obtain the trigonometric Calogero--Moser--Sutherland model of
type $BC_N$, and we can investigate the Inozemtsev model of type
$BC_N$ by perturbation from the trigonometric model
\cite{KT,TakHI}.

A method of quasi-solvability is available on the Inozemtsev model
of type $BC_N$. Finkel et al. studied quasi-solvable models in
\cite{FGGRZ1,FGGRZ2}, and they found several quasi-exactly
solvable many-body systems including the Inozemtsev model of type
$BC_N$. We now describe the f\/inite-dimensional spaces which are
related to the quasi-solvability. The quasi-solvability with
respect to the Hamiltonian $H$ was established in \cite{FGGRZ2}
and reformulated in \cite{Takq}.
\begin{proposition}[\cite{FGGRZ2, Takq}] \label{prop:Hinv}
Let $a$, $b_i$ $(i=0,1,2,3)$ be the numbers which satisfy $a\in \{
-l, l+1 \}$ and $b_i \in \{-l_i/2, (l_i+1)/2 \}$ $(i=0,1,2,3)$.
Set
\begin{gather*}
 \Phi(z)=\prod_{1\leq j<k\leq N} (z_j-z_k)^a \prod_{j=1}^N \prod _{i=1}^3 (z_j-e_i)^{b_i}.  
\end{gather*}
Assume that $d=-((N-1)a+b_0+b_1+b_2+b_3)$ is a non-negative
integer. Let $W^{\rm sym}_d$ be the space spanned by
\begin{gather*}
\Phi(\wp (x_1), \wp (x_2), \dots ,\wp (x_N))\sum _{\sigma \in S_N}
\wp (x_1)^{m_{\sigma (1)}} \wp (x_2)^{m_{\sigma (2)}} \cdots \wp
(x_N)^{m_{\sigma (N)}}
\end{gather*}
such that $m_i \in \{ 0,1, \dots ,d\}$ for all $i$. Then we have
\begin{gather*}
H  \cdot W_d^{\rm sym} \subset W_d^{\rm  sym}.
\end{gather*}
\end{proposition}
The quasi-solvability was extended to the commuting dif\/ferential
operators.
\begin{proposition}[\cite{Takq}, Theorem 3.3] \label{thm:Pinv}
Assume that $d=-((N-1)a+b_0+b_1+b_2+b_3)$ is a non-negative
integer. Then $P_k  \cdot W_d^{\rm sym} \subset W_d^{\rm  sym}$
for $k=1,2,\dots ,N$, where $W_d^{\rm sym}$ is the
finite-dimensional space defined in Proposition~{\rm
\ref{prop:Hinv}} and $P_k$ are the commuting differential
operators which ensure the complete integrability.
\end{proposition}
By the quasi-solvability, f\/initely-many eigenvalues and
eigenfunctions are calculated by diagonalizing the commuting
matrices simultaneously for the case that the assumption of
Proposition~\ref{prop:Hinv} is true. The eigenfunctions obtained
by the quasi-solvability may not be square-integrable in general,
although the eigenfunctions for the case that the parameters $a$,
$b_0$, $b_1$ in Proposition~\ref{prop:Hinv} are chosen as $a=l+1$,
$b_0=(l_0+1)/2$ and $b_1=(l_1+1)/2$ are square-integrable.

It seems that an explicit expression of the Bethe Ansatz as
Proposition \ref{prop:BAAN-1} for the Inozemtsev model of type
$BC_N$ and corresponding conformal f\/ield theory are not known in
the moment of writing, although Chalykh, Etingof and Oblomkov
\cite{CEO} gave a general recipe for calculating the Bloch
eigenfunctions. They showed that these are parametrized by a
certain algebraic variety (the Hermite--Bloch variety) which can
be computed. This would lead to a version of the Bethe Ansatz for
the models including the Inozemtsev model of type $BC_N$, though
these Bethe Ansatz equations would be rather complicated. For a
special case of the $BC_2$ case, this scheme is worked out
explicitly in \cite[\S~6]{CEO}.

We hope to investigate the Bethe Ansatz for the Inozemtsev model
of type $BC_N$ to study the model in more detail.

\section[Towards finite-gap integration of the Inozemtsev model]{Towards f\/inite-gap integration of the Inozemtsev model} \label{sec:fingapICMS}

In Section~\ref{sec:Heun}, we reviewed the f\/inite-gap
integration of Heun's equation, and observed that the existence of
commuting operator of odd order plays important roles.

For a multidimensional generalization of f\/inite-gap integration,
Chalykh and Veselov introduced the notion ``algebraic
integrability''. The Schr\"odinger operator $L=
-\sum\limits_{i=1}^N \partial ^2 /\partial x_i^2 +u(x_1 ,\dots ,
x_N)$ is called completely integrable \cite{CV2}, if there exist
$N$ commuting operators $L_1=L ,L_2, \dots ,L_N$ with
algebraically independent constant highest symbols $s_1 (\xi )$
$(=\xi_1^2 +\cdots +\xi_N^2 )$, $s_2 (\xi), \dots ,s_N(\xi)$ $(\xi
_j =\sqrt{-1} \partial /\partial x_j)$. For example, the
Calogero--Moser--Sutherland models are completely integrable. On
the model of type $A_N$, the highest symbols of $P_k$ (see
equation~(\ref{eq:AN-1Pk})) are written as
$((-\sqrt{-1})^k/k!)\sum\limits_{i_1< i_2 < \cdots <i_k} \xi
_{i_1} \xi _{i_2} \cdots \xi _{i_k} $. The Inozemtsev model of
type $BC_N$ is also completely integrable. The operator $L$ is
called algebraically integrable in the sense of~\cite{CV2}, if $L$
is completely integrable and there exists one more operator $L_0$
commuting with $L_i$ $(i=1,\dots ,N)$ and the highest symbol $s_0
(\xi )$ of $L_0$ takes the distinct values at the roots of the
equations $s_i(\xi )=E_i$ $(i=1,\dots ,N)$ for almost all $E_i$.

On Heun's equation in the elliptic form, if $l_0$, $l_1$, $l_2$,
$l_3$ are integers, then it is algebraically integrable, because
there exists a commuting operator $A$ of odd order. Thus we expect
that algebraically integrable Schr\"odinger operators also have
rich properties.

Chalykh and Veselov conjectured \cite{CV1} that the
Calogero--Moser--Sutherland model with integral coupling constants
are algebraically integrable. For the case of type $A_N$,
Braverman, Etingof and Gaitsgory \cite{BEG} obtained algebraic
integrability. More precisely, they established that, if $l$ is a
positive integer, then the operator
\begin{gather*}
H =- \frac{1}{2} \sum_{i=1}^{N} \frac{\partial ^{2}}{\partial
x_{i}^{2}} + l (l+1) \sum_{1 \leq i<j \leq N} \wp ( x_{i}-x_{j})
\end{gather*}
is algebraically integrable by applying the Bethe Ansatz due to
Felder--Varchenko (see Section~\ref{sec:resultAN-1}) and the
dif\/ferential Galois theory. In \cite{CEO}, Chalykh, Etingof and
Oblomkov proved that the Chalykh--Veselov conjecture is true and
the Inozemtsev model of type $BC_N$ (see equation~(\ref{InoHam}))
is also algebraically integrable, if $l$, $l_0$, $l_1$, $l_2$,
$l_3$ are all integers. Their method relies on results on the
dif\/ferential Galois theory obtained in \cite{BEG} and the local
triviality of the monodromy. On an application of the algebraic
integrability, the eigenfunctions of the Baker--Akhiezer (Bloch)
type are considered. We expect further studies for applications of
the algebraic integrability on the Calogero--Moser--Sutherland
models and the Inozemtsev models.

The explicit expressions of the extra commuting operators were
obtained and investigated by Oblomkov, Khodarinova and  Prikhodsky
\cite{Obl,KP1,KP2} for the case $l=1$ on the
Calogero--Moser--Sutherland model of type $A_3$ and the case
$l=l_0=1$, $l_1=l_2=l_3=0$ on the Inozemtsev model of type $BC_2$.
On the Calogero--Moser--Sutherland model of type $A_3$ with $l=1$,
the Hamiltonian and commuting operators which guarantee the
complete integrability are given as
\begin{gather*}
 H= -(\partial _1 ^2 + \partial _2 ^2 +\partial _3 ^2 )/2+2(\wp _{12} +\wp _{23}+\wp _{31}) ,\\
 P_1= \partial _1 + \partial _2 +\partial _3 , \nonumber \\
 P_3= \partial _1 \partial _2 \partial _3 +2 \wp _{12} \partial _3 +2 \wp _{23} \partial _1 +2 \wp _{31} \partial _2 , \nonumber
\end{gather*}
(see equation~(\ref{eq:AN-1Pk})) where we have used the notations
$\partial _i = \partial /\partial x_i$ and $\wp _{ij} =\wp (x_i -
x_j)$. The additional commuting operators are written as
\begin{gather*}
I_{12} =(\partial _1 -\partial _3)^2 (\partial _2 -\partial _3)^2  -8 \wp _{23} (\partial _1 -\partial _3)^2 -8 \wp _{13}(\partial _2 -\partial _3)^2 \\
\phantom{I_{12} =}{}+4(\wp _{12}-\wp _{13}-\wp _{23} ) (\partial _1 -\partial _3 )(\partial _2 -\partial _3) -2(\wp '_{12} +\wp ' _{13}+6\wp ' _{23}) (\partial _1 -\partial _3 ) \nonumber \\
\phantom{I_{12} =}{} -2 (-\wp '_{12} +6\wp ' _{13}+\wp ' _{23}) (\partial _2 -\partial _3 ) -2\wp '' _{12}  -6\wp '' _{13} -6\wp '' _{23} +4(\wp ^2 _{12} +\wp ^2 _{13} +\wp ^2 _{23}) \nonumber \\
\phantom{I_{12} =}{} +8(\wp _{12} \wp_{13}+\wp _{12} \wp_{23}+7\wp
_{13} \wp_{23}) , \nonumber
\end{gather*}
$I_{23}$ and $I_{31}$, which are written by permuting the indices.
Then any non-symmetric linear combination of them, e.g.,
$L_4=I_{12}+2I_{23}$ would f\/it into the def\/inition of
algebraic integrability (see~\cite{KP1}). Explicit expressions of
the extra commuting operators for the models which have symmetry
of the deformed root system $A_3(m)$ or $B_2(l,m)$ were also
obtained.

Another possible method for constructing extra commuting operators
is the multidimensional Darboux transformation, because the
commuting operator for the case of Heun's equation is constructed
by composing the (generalized) Darboux transformations.
Multidimensional Darboux transformations were studied from several
viewpoints \cite{ABI,GK,Sab,Cha}. In \cite{ABI}, based on the
existence of an explicit eigenfunction of the Hamiltonian
$H(=H^{(0)})$ with a certain eigenvalue, an alternate Hamiltonian
$\tilde{H} (=H^{(N)})$, matrix valued operators $H^{(i)}$
$(i=1,\dots ,N-1)$ and supersymmetry operators $Q ^{-} _{j+1,j}$
and $Q ^{+} _{j,j+1}$ $(j=0, \dots ,N-1)$ which connect $H^{(j)}$
and $H^{(j+1)}$ were introduced and studied. On the other hand, we
know an explicit eigenfunction of the Inozemtsev model of type
$BC_N$, if the value $d(=-((N-1)a+b_0+b_1+b_2+b_3))$ $( a\in \{
-l, l+1 \}, b_i \in \{-l_i/2, (l_i+1)/2 \}$ $(i=0,1,2,3))$ in the
assumption of Proposition \ref{prop:Hinv} is equal to zero. Then
the alternate Hamiltonian $\tilde{H}$ with respect to $H$ in
equation~(\ref{InoHam}) would be written as
\begin{gather*}
 \tilde{H} =-\sum_{j=1}^N\frac{\partial ^2}{\partial x_j^2}+2a(a+1)\sum_{j<k} \left( \wp (x_j-x_k) +\wp (x_j +x_k) \right) \\
 \phantom{\tilde{H} =}{}+ \sum_{j=1}^N \sum _{i=0}^3 2b_i(2b_i +1) \wp(x_j +\omega_i).
\end{gather*}
In the moment of writing, we do not know an operator $L$ which
directly intertwines the operators~$H$ and $\tilde {H}$ as
$\tilde{H}L=LH$. We expect applications of the multidimensional
Darboux transformation for the analysis of the elliptic
Calogero--Moser--Sutherland model or the Inozemtsev model.

\subsection*{Acknowledgements}

The author would like to thank the referees for valuable comments.

\pdfbookmark[1]{References}{ref}
\LastPageEnding

\end{document}